\documentclass[12pt,a4paper]{amsart}

\usepackage{amscd, amssymb,amsmath}

\begin{document}

\pagestyle{headings}

\font\fiverm=cmr5 % Now needed to avoid an error when inputing pictex
\input{prepictex}
\input{pictex}
\input{postpictex}

\def\clsp{\overline{\operatorname{span}}}

\newtheorem{thm}{Theorem}[section]
\newtheorem{cor}[thm]{Corollary}
\newtheorem{lem}[thm]{Lemma}
\newtheorem{prop}[thm]{Proposition}
\newtheorem{thm1}{Theorem}

\theoremstyle{definition}
\newtheorem{dfn}[thm]{Definition}
\newtheorem{dfns}[thm]{Definitions}

\theoremstyle{remark}
\newtheorem{rmk}[thm]{Remark}
\newtheorem{rmks}[thm]{Remarks}
\newtheorem{example}[thm]{Example}
\newtheorem{examples}[thm]{Examples}
\newtheorem{note}[thm]{Note}
\newtheorem{notes}[thm]{Notes}

\numberwithin{equation}{section}

\title{Some Intrinsic Properties of Simple Graph $C^*$-algebras}
\author{David Pask and Seung-Jai Rho}
\address{Mathematics, School of Mathematical and Physical Sciences, The
University of Newcastle, NSW 2308, AUSTRALIA}
\email{davidp@maths.newcastle.edu.au, rho2106@hotmail.com}
\date{\today}

\begin{abstract}
To a directed graph $E$ is associated a $C^*$-algebra $C^* (E)$ called
a graph $C^*$-algebra. There is a canonical action $\gamma$ of ${\bf T}$ on
$C^* (E)$, called the gauge action. In this paper we present
necessary and sufficient conditions for the fixed point algebra $C^*
(E)^\gamma$ to be simple. Our results also yield some structure theorems
for simple graph algebras.
\end{abstract}

\maketitle

\section{Introduction}

\noindent This paper brings together ideas from the theory of
nonnegative matrices associated to strongly connected graphs and
topological graph theory to prove some structure results for graph
$C^*$-algebras. Because of the diverse backgrounds involved, we
have made an effort to make this paper self-contained by including
a little relevant background from each of these areas.

We begin by establishing our notation and conventions for directed
graphs and their $C^*$-algebras (we are aware that graph theorists
use quite different terminology (see \cite{co} for example)). Next
we bring together two sets of results on the simplicity of graph
$C^*$-algebras due to Paterson \cite{p} and Szyma\'{n}ski
\cite{sz}. We show that up to Morita equivalence, a simple graph
$C^*$-algebra is either AF or the $C^*$-algebra of a strongly
connected graph.

We then give some results about the finite path space of a strongly
connected graph. These results are essentially restatements of
standard facts about irreducible nonnegative matrices.
In the following
section we describe the construction of a relative skew product
graph from $E$ (cf.\ \cite{gt,gt1}). Essentially, if $\Gamma$ is a
group with subgroup $H$, a relative skew product graph $E \times_c
( H \backslash \Gamma )$ is an extension of $E$ by the homogeneous
space $H \backslash \Gamma$ using a labelling of the edges in $E$
by elements of $\Gamma$. This construction generalises the usual
skew product graph used in \cite{kp,kqr}. In Theorem
\ref{coverchar} we show that any connected covering graph of a
given connected graph can be written as a relative skew product
graph of the base graph, a generalisation of the results in
\cite{gt,gt1}. We describe an invariant, called the local voltage
group, which enables us to write the connected components of an
ordinary skew product graph as skew product graphs in their own
right. For a row-finite graph $E$ there is a canonical labelling
of the edges in $E$ with the integer $1$ such that $C^* ( E
\times_c {\bf Z} )$ is isomorphic to $C^* ( E) \times_\gamma {\bf
T}$ where $\gamma$ denotes the gauge action of ${\bf T}$ on $C^*
(E)$.

In the final section we apply our results to graph $C^*$-algebras.
In various stages we prove that the fixed point algebra $C^*
(E)^\gamma$ of the gauge action $\gamma$ on $C^* (E)$ is simple if and
only if either $E$ consists of a single vertex or $E$ is row-finite
and has a cofinal subgraph with finitely 
many vertices which is strongly connected with period one. 
If $E$ is strongly connected with
period $d$ then $C^* (E)^\gamma$ is the direct sum of $d$ mutually
isomorphic AF algebras. If, in addition $E^0$ is finite then $C^* (E)$
is stably
isomorphic to a crossed product of a simple AF algebra by a ${\bf Z}$-action.

The authors would like to thank the referee for several helpful
comments
which led to improvements in the final draft.

\section{Graphs and their $C^*$-algebras}

\noindent A {\em directed graph} $E$ consists of a sets $E^0, E^1$ of
vertices 
and edges respectively, together with maps $ r,s : E^1 \rightarrow
E^0$ giving the 
direction of each edge. A {\em subgraph} $F$ of
$E$ consists of subsets $F^i \subseteq E^i$ for $i=0,1$ such that
$s(F^1) \subset F^0$ and $r (F^1) \subset F^0$. A {\em path}
in the directed graph $E$ is a finite sequence $\alpha = \alpha_1
\ldots \alpha_n$ of edges such that $r ( \alpha_i ) = s (
\alpha_{i+1} )$ for $1 \le i \le n-1$. For a finite path $\alpha =
\alpha_1 \ldots \alpha_n$, its {\em length} $\vert \alpha \vert :=
n$ is the number of edges in the sequence $\alpha$. An {\em
infinite path} in $E$ is an infinite sequence $( x_i )_{i \geq 1}$
of edges such that $r ( x_i ) = s ( x_{i+1} )$ for $i \geq 1$; the set
of infinite paths in $E$ is denoted $E^\infty$.
For $n \ge 0$, let $E^n$ denote all those paths in the
directed graph $E$ of length $n$. Let $E^*  = \bigcup_{n \ge 0} E^n$
denote the set of all finite paths in $E$. The
range and source maps extend naturally to $E^*$; for $\alpha =
\alpha_1 \ldots \alpha_n \in E^*$ define $r ( \alpha ) = r (
\alpha_n )$ and $s ( \alpha ) = s ( \alpha_1 )$.
The graph $E$ is called {\em row-finite} if every vertex emits
finitely many edges. A vertex which does not emit any edges is called
a {\em sink}. A vertex which does not receive any edges is called a {\em
source}.

Let $E$ be a directed graph, for $e \in E^1$ we formally denote by
$e^{-1}$, the edge $e$ traversed backwards, so that $s ( e^{-1} )
= r ( e )$ and $r ( e^{-1} ) = s ( e )$. The set of reverse edges
is denoted $E^{-1}$; it is then natural to define $( e^{-1} )^{-1}
=e$ for $e^{-1} \in E^{-1}$.
A {\em walk} in the directed graph $E$ is a sequence  $a = a_1 \cdots a_n $,
where $a_i \in E^1 \cup E^{-1}$ are such
that $r ( a_i ) = s ( a_{i+1} )$ for $i=1 , \ldots , n-1$; we write $s(a) =
s ( a_1 )$ and $r ( a ) = r ( a_n )$.
A walk $a = a_1 \cdots a_n$ is said to be {\em reduced} if it does
not contain the subword $a_i a_{i+1} = a  a^{-1}$ for any $a \in
E^1 \cup E^{-1}$. Given a reduced walk $a = a_1 \cdots a_n$ the
reverse walk is 
written $a^{-1} :=  a_n^{-1} \cdots a_1^{-1}$, which is also
reduced. If $a, b$ are reduced walks with $r(a) = s(b)$, then $a \cdot b$
will be understood to be the reduced walk obtained by concatenation
and then cancellation using
the relations $e e^{-1} = s(e)$, $e^{-1} e =  r(e)$, $e^{-1} s(e) =
e^{-1} = r(e)e^{-1}$ and $e r(e) =
e = s(e) e$ for $e \in E^1$.
With composition and inverse operations defined above, the set $\pi_1
(E)$ of reduced walks, forms a groupoid with unit space $E^0$  and is
referred to as the {\em fundamental groupoid} of $E$ (note that
the roles of the range and source map must be reversed to make
$\pi_1 (E)$ a category, cf.\ \cite{kp2}).

Let $E$, $F$ be directed graphs. A {\em graph morphism} $\phi : F
\rightarrow E$
consists of maps $\phi^i : F^i \rightarrow E^i$ for $i=0,1$ such that
$\phi^0 ( r ( f ) ) = r ( \phi^1 (f) )$ and $\phi^0 ( s ( f ) )
= s ( \phi^1 (f) )$ for all $f \in E^1$.

\begin{dfn}
Let $E$, $F$ be directed graphs and $p : F \rightarrow E$ be a graph
morphism; then $p$ is {\em covering map} if for each $v \in F^0$,  $p$
maps $r^{-1} (v)$ bijectively onto $r^{-1} (p(v))$ and $s^{-1} (v)$
bijectively onto $s^{-1} (p(v))$ (cf.\ \cite[\S 2.1.6]{st}).
\end{dfn}

\noindent Let $p : F \rightarrow E$ be a covering map, for $a =
a_1 \cdots a_n \in \pi_1 (F)$ let $p(a) = p ( a_1 ) \cdots p ( a_n
)$ where $p ( f^{-1} ) := p (f)^{-1}$ for $f \in E^1$. In this
way, the covering map $p$ induces a morphism $p_* : \pi_1(F)
\rightarrow \pi_1(E)$. A graph morphism $p : F \rightarrow E$ has
the {\em unique walk lifting property} (see \cite[Theorem
2.1.1]{gt}) if given $u \in F^0$ and $a \in \pi_1(E)$  with $s ( a
) = p^0 ( u )$, there is a unique $\tilde{a} \in \pi_1 (F)$ such
that $s ( \tilde{a} ) = u$ and $p ( \tilde{a} ) = a$; one may also
lift reduced walks which end at $p^0 (u)$. Likewise one may lift
infinite paths from $E$ to $F$.
The following result is routine (see % \cite[\S 2.2.2]{st},
\cite[Lemma 17.4]{ko} for instance):

\begin{lem} \label{cmapiffupl}
Let $E$, $F$ be directed graphs and $p : F \rightarrow E$ be a graph
morphism. Then $p$ has the unique walk lifting property if and only if
$p$ is a covering map.
\end{lem}

\noindent
The directed graph $E$ is said to be {\em connected} if,
given  any two distinct vertices in $E$, there is a reduced
walk between them. A directed graph $T$ is a {\em tree} if and
only if there is precisely one reduced walk between any two
vertices (so a tree is connected). For a connected graph $E$ and $v \in
E^0$ we may define
\[
\pi_1 (E,v) = \{ a \in \pi_1 (E) : s(a) = v = r (a) \}
\]

\noindent
so that $\pi_1 (E,v)$ is the isotropy group of the unit $v$ in $\pi_1
(E)$. If $E$ is connected, then the groupoid $\pi_1 (E)$
is connected and so all its isotropy groups are
isomorphic. Our definition of the fundamental group matches the usual
one given in \cite[\S 2.1.6]{st} because taking the reduction of a
walk coincides with the notion of path equivalence used there.

Using the axiom of choice it may be shown that every connected graph
$E$ contains a spanning tree $T$ (cf.\ \cite[\S 2.1.5]{st}): a
subgraph $T$ which is itself a tree with $T^0 = E^0$. Fix a
given spanning tree $T$ of $E$ and a vertex $v \in E^0$, then for $w \in
T^0$ we set $b_w$ to be the unique reduced walk in $T$ from $v$
to $w$. If $E$ has a spanning tree then it is connected.

Let $E$ be a directed graph, then a {\em Cuntz-Krieger $E$-family}
(or a {\em representation} of $E$) consists of a family $\{ P_v :
v \in E^0 \}$ of mutually orthogonal projections and a family $\{
S_e : e \in E^1 \}$ of partial isometries with mutually orthogonal
ranges such that
\[
S_e^* S_e = P_{r(e)} , \quad
S_e S_e^* \le P_{s(e)} , \quad
P_v = \sum_{s(e)=v} S_e S_e^* \text{ if } 0 < |s^{-1} (v) | < \infty.
\]

\noindent The graph $C^*$-algebra $C^* (E)$ is generated by a
universal Cuntz-Krieger $E$-family $\{ s_e , p_v \}$ (see \cite[\S
1]{rs}). The class of
graph $C^*$-algebras is quite broad, it includes, up to Morita
equivalence, all AF algebras and all Cuntz-Krieger algebras (see
\cite{kpr}, \cite{kprr} amongst others).

\section{Simple graph algebras}

\noindent A {\em loop} in $E$ is a path $\alpha$ with $\vert
\alpha \vert \ge 1$ such that $s ( \alpha ) = r ( \alpha  )$.
% in graph theory these are called circuits;
% moreover a loop is a circuit which consists of a single edge.
The loop $\alpha \in E^n$ is {\em simple} if the vertices
$\{ r ( \alpha_i ) : 1 \le i \le n \}$ are distinct.
% in graph theory such a circuit is called a cycle.
The graph $E$ satisfies condition (K) if no
vertex lies on exactly one simple loop.
The directed graph $E$ is {\em cofinal} if, given $x
\in E^\infty$ and $v \in E^0$ there is a path $\alpha$ such that
$s ( \alpha )= v$ and $r ( \alpha ) = r ( x_n )$
for some $n \ge 1$.

If $v , w \in E^0$ then we write $v \ge w$ if there is a path from
$v$ to $w$. A subset $H$ of $E^0$ is {\em hereditary} if $v \in H$
and $v \ge w$ implies that $w \in H$. A hereditary subset $H$ is
{\em saturated} if there is no vertex $v \in E^0 \backslash H$
with $0 < |s^{-1} (v) | < \infty$ such that $r (e) \in H$ for all
$e \in s^{-1} (v)$. For $X \subseteq E^0$ define $L_X$ to be the
hereditary subset consisting of all vertices which $X$ connects
to; then $\Sigma ( X )$, the smallest saturated hereditary subset
containing $X$ consists of all those vertices of $v \in E^0$ with
$0 < |s^{-1} (v) | < \infty$ such that every infinite path
starting at $v$ eventually uses vertices in $L_X$. Recall from
\cite{bhrs} that a nontrivial saturated hereditary subset $H$ of
$E$ gives rise to a nontrivial gauge invariant ideal $I_H$ of $C^*
(E)$.

\begin{thm} \label{ps}
Let $E$ be a directed graph then $C^* (E)$ is simple if and only if
\begin{itemize}
\item[(i)] $E$ is cofinal,
\item[(ii)] $E$ satisfies condition (K),
\item[(iii)] If vertex $v$ emits infinitely many edges then every vertex
connects to $v$.
\end{itemize}
\end{thm}

\begin{proof}
See \cite[Theorem 4]{p} or \cite[Theorem 12]{sz}. In \cite[Theorem 4]{p} the
proof of the only if direction is omitted, but as we see below it is
not hard. 

Suppose that $E$ is not cofinal, then there are $x \in E^\infty$ and
$v \in E^0$ such that $v$ does not connect to any vertex used by $x$. Hence
$s(x) \not\in \Sigma ( L_{\{ v \}} )$ because the infinite path $x$ does not
enter $L_{\{v\}}$. Therefore $\Sigma ( L_{\{v\}} )$ gives rise to
a nontrivial gauge invariant ideal $I_{\Sigma ( L_{\{v\}} )}$ of $C^*
(E)$ in which case $C^* (E)$ is not simple.

Suppose that $E$ does not satisfy condition (K), then there is a
simple loop $L$  in $E$ such that each vertex of $L$ lies on no
other loop. Either $E$ is cofinal or it is not. If $E$ is not
cofinal then $C^* (E)$ is not simple by the argument above. If $E$
is cofinal then by \cite[Theorem 2.4]{kpr} $C^* (E)$ is strongly
Morita equivalent to $C ( {\bf T} ) $ which is not simple.

If $E$ has a vertex $v$ of infinite valency and a vertex $w$ which does not
connect to $v$, then $v \not\in L_{\{w\}}$; moreover $v \not\in \Sigma
( L_{\{w\}} )$,
since $v$ has infinite valency. Hence $\Sigma ( L_{\{w\}} )$ is
a nontrivial saturated hereditary subset of $E^0$ and so gives rise to
a nontrivial gauge invariant ideal $I_{\Sigma ( L_{\{w\}} )}$ of $C^*
(E)$. Hence $C^* (E)$ is not simple, which concludes the proof.
\end{proof}

\noindent The directed graph $E$ is {\em strongly connected} if
for every pair of vertices $v, w$ there is a path $\alpha$ with $\vert
\alpha \vert \ge 1$
such that $s ( \alpha ) = v$, $r ( \alpha ) = w$. A
strongly connected graph is sometimes said to be transitive or
irreducible (since its vertex incidence matrix $A_E$ is irreducible
(cf. \cite[\S 1.3]{s})). If $E$ is strongly connected, then it is
cofinal. A subgraph $F$ of $E$ is said to be {\em
cofinal} if for every $x \in E^\infty$ there exists $N(x)$ such that
$r ( x_n ) \in F^0$  for $n \ge N(x)$.

\begin{thm} \label{reducetosc}
If the directed graph $E$ is cofinal then either $E$ has no loops, or
there exists a strongly connected cofinal subgraph $F \subseteq E$.
\end{thm}

\begin{proof}
Under the relation of mutual connectivity (i.e.\ $v \sim w$ if and
only if $v \ge w$ and conversely), there exists an equivalence
class $X$ of vertices which contains any vertex which lies on a
loop. If $E$ has a loop and is cofinal there can only be one such
equivalence class. Let $F$ be the subgraph with $F^0 = X$ and $F^1
= \{ e : s(e) \in F^0 \}$. We claim that $r ( F^1) \subseteq F^0$.
Suppose that $e \in F^1$ is such that  $r(e) \not\in F^0$. Since
$s(e) \in F^0$, the vertex $s(e)$ connects to itself via $\alpha$.
Let $y = \alpha \alpha \alpha \cdots \in E^\infty$ then since $E$
is cofinal $r(e)$ must connect to some vertex on $\alpha$ and
hence to $s(e)$. But this contradicts our assumption that $r(e)
\not\in F^0$. By its definition, $F$ is strongly connected, it
remains to show that $F \subseteq E$ is cofinal.

Let $x \in E^\infty$ then if $s ( x ) \in F^0$ then we may take
$N(x)=1$. So we suppose that  $s(x) \not\in F^0$. Let $v \in F^0$,
then there is a loop $\alpha$ with $s( \alpha )  = v$ which yields
$y = \alpha \alpha \cdots \in E^\infty$. Since $E$ is cofinal there
exists $\beta \in E^*$ such that $s ( \beta ) = v$ and $r ( \beta ) =
r ( x_m )$ some $m \geq 1$. Let $w = r ( x_{m+1} )$ then
by cofinality there must be a path $\gamma \in E^*$ with $s ( \gamma ) = w$
and $r ( \gamma ) = r ( y_{\ell \vert \alpha \vert} ) = v$ for some
$\ell \ge 1$. Then $v$ connects to $w$ via $\beta x_{m+1}$ and $w$
connects to $v$ via $\gamma$, and so $w \in F^0$ by definition.
\end{proof}

\begin{cor} \label{ga-appl}
Suppose that $C^* ( E )$ is simple then either $C^* ( E )$ is an
AF-algebra or there exists a strongly connected subgraph $F \subseteq
E$ such that $C^* ( E )$ is strongly Morita equivalent to $C^* ( F )$.
\end{cor}

\begin{proof}
If $C^* (E)$ is simple, then by Theorem \ref{ps} the graph $E$ is
cofinal. If $E$ has no loops then by Theorem \ref{ps} $E$ cannot have
any vertices of 
infinite valency and then $C^* (E)$ is an AF-algebra by \cite[Theorem
2.4]{kpr}. On the other hand if $E$ has loops,
then by Theorem \ref{reducetosc} there exists a strongly connected
cofinal subgraph $F$ of $E$. Let $P = \sum_{v \in F^0}
p_v$, then by a similar argument to the one given in \cite[\S 1]{bprs}
$P$ defines a projection in the multiplier algebra of $C^*
(E)$ such that $P C^* (E) P \cong C^*  (F)$. To see that $P C^* (E) P$
is full, just observe that $P \neq 0$ and that $C^* (E)$ is simple. Hence
$C^* (E)$ is strongly Morita equivalent to $C^* (F)$ as required.
\end{proof}

\section{Strongly connected graphs} \label{nnm}

\noindent In this section we briefly give some results about the
path space of a strongly connected graph. Essentially, they are
reformulations of results about finite, nonnegative matrices.
% -- in graph theory such graphs are said to have infinite order.
For more details about nonnegative matrices see \cite[\S 1.3]{s}
or \cite[\S 4]{lm}. Let $E$ be a directed graph, then for $v \in
E^0$ we define its {\em period} $d(v)$, to be the greatest common
divisor of the lengths of all loops which begin at $v$. If there
are no such loops then we set $d(v)=0$.

\begin{lem} \label{dv=dw}
If $E$ is strongly connected then for any $v, w \in E^0$ we have $d (
v)  = d ( w )$.
\end{lem}

\begin{proof}
Let $v, w \in E^0$, then since $E$ is strongly connected
there are paths $\alpha \in E^k$ and $\beta \in E^\ell$ from $v$ to
$w$ and $w$ to $v$, respectively. So, if $w$ is the source of a loop of
length $s$ then  $v$ is the source of loops of length $k + \ell + s$
and $k + \ell + 2s$. Hence $d(v)$ divides $( k + \ell + 2s ) - ( k +
\ell + s ) = s$.
Therefore $d(v)$ divides $s$ for every loop with source $w$ of length
$s$. Since $d(w)$ is the greatest common divisor of such numbers $s$ we must
have $d(v) \le d(w)$. But since the argument can be repeated with $v$ and $w$
exchanged we have $d(w) \le d(v)$ and so $d(v)=d(w)$ as required.
\end{proof}

\noindent
Hence we may define the {\em period} of a strongly directed graph to
be the period of any one of its vertices. If $E$ is strongly connected
with period $1$ and has finitely many vertices then  its vertex
incidence matrix $A_E$ is aperiodic in the sense
that there is a $k \ge 1$ such that every entry of $A_E^k$ is strictly
positive (cf.\ \cite[p.253]{ck}).

\begin{lem} \label{eventual}
Let $E$ be a strongly connected graph with period $d$. For each
$v \in E^0$ there is a positive integer $N (v)$ such that for $k \geq
N (v)$, $v$ is the source of a loop of length $kd$.
\end{lem}

\begin{proof}
Suppose that $\alpha \in E^{kd}$ and $\beta \in E^{sd}$ are loops with
source $v$. Then $\alpha \beta \in E^{( k + s) d}$ is a loop with
source $v$. Hence the set
\[
V = \{ kd  : \text{ there is a loop of the length $kd$ in $E$ with
source $v$ } \}
\]

\noindent
of positive integers is closed under addition, and also their the
greatest common divisor
is $d$. But then by Lemma \cite[Lemma A3]{s} $V$ must contain all
but a finite number of positive multiples of $d$, and the result follows.
\end{proof}

\noindent
Next we examine the structure of the finite path space of a strongly
connected graph with period $d$ (see \cite[Theorem 1.3 in Part I]{s}):

\begin{lem} \label{resclass}
Let $E$ be a strongly connected graph with period $d$, and $v
\in E^0$. For any $w \in E^0$ there exists
 $r_w$ with $0 \leq {r_w} < d$ such that
\begin{itemize}
\item[(i)] if $\mu \in E^s$ is a path from $v$ to $w$ then $s \equiv
{r_w}$ (mod $d$),
\item[(ii)] there is a positive integer $N(w)$ such
that for $k \ge N(w)$ there is a path of length $kd+{r_w}$ from $v$ to $w$.
\end{itemize}
\end{lem}

\begin{proof}
For (i), let $\alpha \in E^m$ and
$\beta \in E^{m'}$ be paths from $v$ to $w$. Since $E$ is
strongly connected there is a path $\gamma \in E^p$ from
$w$ to $v$. Hence $\alpha \gamma \in E^{m+p}$ and $\beta \gamma \in
E^{m'+p}$ are loops with source $v$. Since $E$ has period $d$ it
follows that $d$ divides $m+p$
and $m' +  p$, and so $d$ divides $(m+p) - (m'+p) = m-m'$,
that is, $m - m' \equiv 0$ (mod $d$). Set $r_w = m$ (mod $d$). If
$\mu \in E^s$ is a path from $v$ to $w$ then $s \equiv r_w$ (mod $d$).

For (ii), let $\alpha \in E^n$ be a path from $v$ to $w$, then from
(i) there is a  positive integer $m$ such that $n = md + r_w$.
By Lemma \ref{eventual} there is a positive integer $N_0$ such that
for all $s \ge N_0$, $v$ is the source of a loop of length $sd$.
Set $N(w) = N_0 +m$, if $k \geq N(w)$ then $k-m \ge N_0$ and so
there is a loop $\delta \in E^{(k -m )d}$ with source $v$. Hence
$\alpha \delta$ is a path of length
$n + (k-m)d = (m+k-m) d + r_w = kd +r_w$
from $v$ to $w$.
\end{proof}

\noindent
Note that if $w=v$ in the above Lemma, then $r_v =0$. Note also that $r_w$ in
Lemma \ref{resclass} is called a {\em residue class} of $w$ with
respect to $v$ (cf.\ \cite[Definition 1.7 in Part I]{s}).

\section{Skew product graphs}

\noindent In this section we describe the relative skew product
graph construction and show how it can be used to describe all
connected coverings of a given connected graph. We introduce an
invariant, called the local voltage group which enables us to
describe the connected components of ordinary skew product graphs
and exhibit them as ordinary skew product graphs in their own
right. Finally we specialise to a certain skew product graphs
formed from the integers with a view to applications to graph
$C^*$-algebras in the next section.

The following definition generalises the ones given in
\cite[\S 4]{gt1} and \cite[\S 2.3.2]{gt} for row-finite graphs
with finitely many vertices and groups which are finite. Let $E$ be a
directed graph, $\Gamma$ a countable group with subgroup $H$ and
$c : E^1 \rightarrow \Gamma$ a function, which we think of as a labelling
of the edges in $E$ by elements of $\Gamma$.  From this data we may
form the {\em relative skew product graph} $E \times_c ( H
\backslash \Gamma )$ has vertex set $E^0 \times ( H \backslash
\Gamma )$, edge set $E^1 \times ( H \backslash \Gamma)$ where the
range and source of the the edge $(e,Hg)$ are given by
\[
s ( e, Hg ) = ( s(e), Hg )  \text{ and }  r( e, Hg ) = (r(e), Hg
c(e) ) \text{ respectively}.
\]

\noindent
If $H$ is the trivial subgroup of $\Gamma$ then the above definition
reduces to that of the ordinary skew product graph $E \times_c \Gamma$
used in \cite[Definition 2.1]{kp} and \cite[\S 2.1.1]{gt} (as opposed
to those used in \cite{kqr,dpr} where the labelling gives rise to a coaction
of $\Gamma$ on $C^* (E)$). If $H$ is trivial there is a free action
$\lambda$ of $\Gamma$ on $E \times_c \Gamma$ defined for $g, h \in \Gamma$
by
$\lambda^i_h ( x, g ) = (x, hg ) \text{ where } x \in E^i \text{
for } i = 0 ,1$.
If $H$ is a subgroup of $\Gamma$, then $H$ also acts freely on $E
\times_c \Gamma$ and the quotient graph $H \backslash ( E \times_c
\Gamma )$ is isomorphic to $E \times_c ( H \backslash \Gamma )$
via the maps $[ x , g ] \mapsto (x , Hg )$ where $x \in E^i$ for
$i=0,1$.

For functions $c_1 , c_2 : E^1 \rightarrow \Gamma$ we say that
$c_1$ and $c_2$ are {\em cohomologous} and write $c_1 \sim c_2$ if
there is a function $b : E^0 \rightarrow \Gamma$ such that $c_1
(e) b (r(e)) =  b (s(e)) c_2 (e)$ for all $e \in E^1$. The
relation $\sim$ amongst all functions from $E^1$ to $\Gamma$ is an
equivalence relation. If $c_1 ,  c_2 : E^1 \rightarrow \Gamma$ are
cohomologous then $E \times_{c_1}
\Gamma$ is equivariantly isomorphic to $E \times_{c_2} \Gamma$ and
$E \times_{c_1} ( H \backslash \Gamma ) \cong E \times_{c_2} ( H
\backslash  \Gamma )$ via the maps $(v,g) \mapsto (v, g b(v) )$
and $(e,g) \mapsto (e, g b (s(e)) )$.

The map $p_c : E \times_c ( H \backslash \Gamma) \rightarrow E$
defined by
$p_c ( x, Hg )= x \text{ where } x \in E^i \text{ for } i = 0 ,1$
is a covering map. If $E , E \times_c ( H \backslash
\Gamma )$ are connected and $H$ is normal in $\Gamma$ then $p_c$
is a {\em regular} covering since the image of $\pi_1 ( E \times_c
( H \backslash \Gamma ) , (v, Hg) )$ under $( p_c )_*$
is a normal subgroup of $\pi_1 (E,v)$. In fact, we shall show that
every connected covering of a given connected graph $E$ is a
relative skew product graph formed from $E$ (cf.\ \cite{gt,gt1}).

\begin{lem} \label{emptytree}
Let $E$ and $F$ be connected graphs, $p : F \rightarrow E$ a covering
map, $T$ a spanning tree for $E$ and $v \in E^0$. Then there is a
spanning forest $\{ \tilde{T}_u : u \in p^{-1} (v) \}$ of $F$ such
that $p ( \tilde{T}_u ) = T$ % $u \in \tilde{T}_u$
for all $u \in p^{-1} (v)$.
\end{lem}

\begin{proof}
For $w \in E^0$ let $b_w$ be the unique reduced walk in $T$ from $v$ to
$w$. For $u \in p^{-1} (v)$ let $\tilde{b}_w (u)$ be the unique lift
of $b_w$ to a reduced walk beginning at $u \in F^0$. Put $\tilde{T}_u^0 = r (
\tilde{b}_w (u) )$ and let $\tilde{T}_u^1$ consist of all the edges
comprising each $\tilde{b_w} (u)$ as $w$ runs through $E^0$. By the
unique walk lifting property $\tilde{T}_u$ is  a tree such that $p (
\tilde{T}_u ) = T$.

To show that $\{ \tilde{T}_u : u \in p^{-1} (v) \}$ is a forest we
must show that the $\tilde{T}_u$ are disjoint. Suppose that $w$ is a
vertex in $\tilde{T}_u \cap \tilde{T}_{u'}$, where $u , u' \in p^{-1}
(v)$. Let $a$ be the unique reduced walk in $\tilde{T}_u$ from $u$ to $w$
and $b$ be the unique reduced walk in $\tilde{T}_{u'}$ from
$u'$ to $w$. Then $a b^{-1}$ is a reduced walk in
$F$ from $u$ to $u'$, and so $p (a b^{-1} ) = p(a) p(b)^{-1}$ is a loop in
$T$ with source $v$. Hence $p(a)=p(b)$ and then $a^{-1}=b^{-1}$ by unique
walk lifting at $w$. Therefore we have that $u = s (a) = s (b) =
u'$. A similar argument shows that $\tilde{T}_u$ and $\tilde{T}_{u'}$
cannot have an edge in common.

To show that $\{ \tilde{T}_u : u \in p^{-1} (v) \}$ spans $F$, suppose
that $w \in F^0$, then there exists a unique reduced walk $b_{p(w)}$ in $T$
from $v$ to $p (w )$. By the unique walk lifting property there is a
reduced walk $\tilde{b}_{p(w)}^{-1}$ in $F$ such that
$s ( \tilde{b}_{p(w)}^{-1} ) = w$ and $p ( \tilde{b}_{p(w)}^{-1} ) =
b_{p(w)}^{-1}$. Since $p ( r ( \tilde{b}_{p(w)}^{-1} ) ) = v$ we must
have $r ( \tilde{b}_w^{-1} ) \in p^{-1} (v)$ and so $w$ lies in
$\tilde{T}_{r ( \tilde{b}_w^{-1} )}$.
\end{proof}

\begin{thm} \label{coverchar}
Let $E, F$ be connected graphs, $p : F \rightarrow E$ a covering
map, $T$ a spanning tree for $E$ and $v \in E^0$ . Let $w \in F^0$
be such that $p(w)=v \in E^0$ then there exists a map $c = c_{v,T}
: E^1 \rightarrow \pi_1 (E,v)$ such that $F \cong E \times_{c} (
p_* \pi_1 (F,w) \backslash \pi_1 (E,v) )$.
\end{thm}

\begin{proof}
Let $H = p_* \pi^1 ( F , w )$ and $\{ a_i : i \in I \} \subset
\pi^1 ( E , v )$ be a right transversal for $H$ in $\pi^1 ( E, v
)$. For $i \in I$, let $\tilde{a_i}$ be the lift of $a_i$ in $F$
with $s (\tilde{a}_i ) = w$, put $\tilde{v}_i = r ( \tilde{a}_i )$
then  $p^{-1} (v) = \{ \tilde{v}_i  : i \in I \}$.  Let $c =
c_{v,T} : E^1 \rightarrow \pi^1 (E,v)$ be defined by $c (e) =
b_{s(e)} e b_{r(e)}^{-1}$ where $b_{s(e)}$ and $b_{r(e)}$ are the
unique reduced walks in $T$ from $v$ to $s(e)$ and $r(e)$
respectively.

For $i \in I$ let $\tilde{T}_i$ be the tree $\tilde{T}_{\tilde{v}_i}$
in $F$ described in Lemma \ref{emptytree}. Define a map $\phi : F
\rightarrow E \times_{c}  ( H \backslash \pi^1 ( E , v ) )$ as follows: if
$u \in \tilde{T}_i^0$, then set $\phi^0 ( u ) = (  p ( u  ) , H a_i
)$. If $f \in F^1$ is such that $s ( f ) \in \tilde{T}_j^0$ and $r (
f ) \in \tilde{T}_i^0$, then we put $\phi^1
(f ) = (  p ( f ) ,  H a_i c ( p ( f ) )^{-1}  )$.

For $f \in F^1$
we have $\phi^0 ( s ( f ) ) = ( p ( s ( f ) ) , H a_j )$ whereas
\[
s ( \phi^1 ( f ) ) = s (  p (f ) , H a_i c ( p ( f ) )^{-1}   ) =
(  s ( p ( f )  ) , H a_i c ( p ( f ) )^{-1}   ) ,
\]

\noindent so $s ( \phi^1 (f) ) = \phi^0 ( s(f) )$ provided that $H
a_i c ( p ( f ) )^{-1} = H a_j$. To see this, let
\[
h =  \tilde{a}_i \tilde{b}_{r(f)}^{-1} f^{-1} \tilde{b}_{s(f)}
\tilde{a}_j^{-1}
\]

\noindent
where $\tilde{b}_{r(f)}$ is the lift of $b_{p(r(f))}$ with $s (
\tilde{b}_{r(f)} ) = \tilde{v}_i$ and $\tilde{b}_{s(f)}$ is the lift of
$b_{p(s(f))}$ with $s ( \tilde{b}_{r(f)} ) = \tilde{v}_j$; so $h$ is a
closed walk in $F$ with $s ( h ) = w$. Hence
\[
p ( h ) = a_i b_{p(r(f))}^{-1} p(f)^{-1} b_{p(s(f))} a_j^{-1} =
a_i c ( p(f) )^{-1} a_j^{-1} \in H ,
\]

\noindent and so $H a_i c ( p ( f ) )^{-1}  = H a_j$ as required.
It is straightforward to show that $r  \circ \phi^1 = \phi^0 \circ
r$, and hence $\phi$ is a graph morphism.

We now show that $\phi$ is injective.
Suppose that for some $\tilde{e}, \tilde{f} \in F^1$ we have
$\phi^1 ( \tilde{e} ) = \phi^1 ( \tilde{f} )$, where $r (
\tilde{e} ) \in \tilde{T}_i$, $r ( \tilde{f} ) \in
\tilde{T}_{i'}$, $s ( \tilde{e} ) \in \tilde{T}_j$ and $s (
\tilde{f} ) \in \tilde{T}_{j'}$. Then $(  p ( \tilde{e} ) ,  H a_i
c ( p ( \tilde{e} ) )^{-1}  ) = (  p ( \tilde{f} ) , H a_{i'} c (
p ( \tilde{f} ) )^{-1}
 )$ in which case $p (\tilde{e} ) = p ( \tilde{f} )$ and then
$i=i'$. Let $\tilde{b}$ be a reduced walk in $\tilde{T}_i$ from $s
( \tilde{e} )$ to $s ( \tilde{f} )$ then $b = p ( \tilde{b} )$ is
a closed reduced walk in $T$ which implies that $s ( \tilde{e} ) =
s ( \tilde{f} )$. So $\tilde{e} = \tilde{f}$ by unique walk
lifting of $p ( \tilde{e} ) = p ( \tilde{f} )$ at $s ( \tilde{e}
)$. The proof that $\phi$ is injective on $F^0$ follows similarly.

Finally we claim that $\phi$ is surjective.
Given $( e , H a_j )$ in $( E \times_{c} ( H \backslash \pi_1 (
E,v ) ) )^1$, let $\tilde{a} = \tilde{a_j} \tilde{b}_{s(e)}$ be
the lift of $a_j$ at $w$ followed by the lift of $b_{s(e)}$ at
$\tilde{v}_j$. Set $\tilde{e}$ to be the lift of $e$ in $F$ based
at $r ( \tilde{a} )$, then $s ( \tilde{e} ) \in \tilde{T}_j$ and
suppose that $r ( \tilde{e} ) \in \tilde{T}_i$. Then since $p (
\tilde{e} ) = e$ we have $\phi^1 ( \tilde{e} ) = ( e , H a_i
c(e)^{-1} )$ and $H a_i c(e)^{-1} = H a_j$ by a similar argument
to the one given to show that $\phi^0 \circ s = s \circ \phi^1$
above. The proof that $\phi^0$ is surjective follows similarly.
\end{proof}

\noindent
The function $c = c_{v,T} : E^1 \rightarrow \pi_1 (E,v)$ given in the
proof of Theorem \ref{coverchar} depends on the transversal $p_* \pi_1 (F,w)$
in $\pi_1 (E,v)$, and the choice of $v$ and $T$. In Remarks \ref{indep} we
show that the cohomology class of $c$ is independent of the choice of
$v$ and $T$. By unique walk lifting, a different transversal
has no effect other than to permute the trees in the spanning forest of $F$.

We shall now fix our attention on ordinary skew product graphs. Let
$E$ be a directed graph and $c : E^1 \rightarrow \Gamma$ a function
where $\Gamma$ is a countable group. Then there is a map (which we
shall also denote by $c$) from $\pi_1 (E)$ to $\Gamma$ defined by
\[
c ( a_1 \ldots a_n ) = c ( a_1 ) \ldots c ( a_n ) ,
\]

\noindent
where $c ( e^{-1} ) = c(e)^{-1}$. It is straightforward to show that $c$ is a
functor from $\pi_1 (E)$ to $\Gamma$. For $v \in E^0$ let
\[
\Gamma_v (c) = \{ c ( a ) : a \in \pi^1 ( E , v ) \} .
\]

\noindent
Since it is the range of a homomorphism, $\Gamma_v
(c)$ is a subgroup of $\Gamma$, called the {\em local voltage group}
of $c$ based at $v$ (in \cite[\S 2.5]{gt} these groups are
referred to as the local group). The following facts are not difficult to
establish (see \cite[Lemma 6.5]{r}):

\begin{lem} \label{lvprops}
Let $E$ be a directed graph and $c : E^1 \rightarrow \Gamma$ a function where
$\Gamma$ is a countable group.
\begin{trivlist}
\item[(i)] If there is a walk between $v , w \in E^0$ then $\Gamma_v (c)$ and
$\Gamma_w (c)$ are conjugate subgroups of $\Gamma$.
\item[(ii)]  If $c_1 , c_2 : E^1 \rightarrow \Gamma$ are cohomologous
then for all $v \in E^0$, $\Gamma_v (  c_1 )$ and $\Gamma_v ( c_2
)$ are conjugate subgroups of $\Gamma$.
\end{trivlist}
\end{lem}

\begin{rmk} \label{walks}
Let $E$ be a directed graph and $c : E^1 \rightarrow \Gamma$ a function where
$\Gamma$ is a countable group. If $\tilde{a}$ is a reduced walk in $E
\times_c \Gamma$ from $(v,g)$ to $(u,h)$, then $h = g c(a)$ and $a = p_c (
\tilde{a} )$ is a reduced walk in $E$ from $v$ to $u$
(cf. \cite[Theorem 2.1.2]{gt}).
\end{rmk}

\begin{prop} \label{logrcor}
Let $E$ be a connected graph,  $c : E^1 \rightarrow \Gamma$ a function where
$\Gamma$ is a countable group. Then the vertices $(v,g )$ and $(u,h)$
lie in the 
same connected component of $E \times_c \Gamma$ if and only if there
is a reduced walk $a$ in $E$ from $v$ to $u$  such that $g^{-1} h$
lies in the coset $\Gamma_v ( c ) c(a)$ of the local voltage group of
$c$ based at $v$.
\end{prop}

\begin{proof}
If $(v,g)$ and $(u,h)$ lie in the same connected component of $E \times_c
\Gamma$ then there exists a reduced walk $\tilde{a}$  in $E \times_c
\Gamma$ from $(v,g )$ to $(u,h)$. By Remark \ref{walks} we have  $h =
g c(a)$, where $a = p_c ( \tilde{a} )$, then $g^{-1} h  = c (a)$ lies
in $\Gamma_v ( c ) c ( a )$ as required.

Conversely let $a$ be a reduced walk in $E$ from $v$ to $u$ such that
$g^{-1} h \in \Gamma_v ( c ) c ( a )$. Hence there is $b \in \pi^1 (
E, v )$ such that $g^{-1} h  = c (b) c (a) = c (ba)$. Let
$\tilde{d}$ be the unique lift of $d=ba$ in $E \times_c \Gamma$
starting at $(v,g)$ then $\tilde{d}$ is a reduced walk in $E \times_c \Gamma$
with  $s ( \tilde{d} ) =  (v,g)$ and  $r ( \tilde{d} ) = (  r (ba) , g
c(ba) ) = (u,h )$ so that $(v,g)$ and $(u,h)$ are  in the same
connected component.
\end{proof}

\begin{cor} \label{my2-2}
Let $E$ be a connected graph and $c : E^1 \rightarrow \Gamma$ a
function where 
$\Gamma$ is a countable group. Then $E \times_c \Gamma$ consists of  $ \vert
\Gamma : \Gamma_v ( c ) \vert$ mutually isomorphic connected
components, where 
$\Gamma_v (c)$ is the local voltage group of $c$ based at $v \in
E^0$. In particular, $E \times_c \Gamma$ is connected if and only if
$\Gamma_v ( c ) = \Gamma$ for some $v \in E^0$.
\end{cor}

\begin{proof}
By Proposition \ref{logrcor} the number of connected components of $E
\times_c \Gamma$ is equal to  $\vert \Gamma : \Gamma_v ( c ) \vert$,
the number of cosets of $\Gamma_v (c)$ in $\Gamma$. Since the natural action
$\lambda$ of $\Gamma$ on $E \times_c \Gamma$  is transitive on
$p_c^{-1} (v) = 
\{ ( v , g ) : g \in \Gamma \}$ for all $v \in E^0$ and $E$ is
connected it follows that the connected components of $E \times_c
\Gamma$ are mutually isomorphic.
\end{proof}

\noindent
There appears to be an error on \cite[p.88]{gt} where it is shown that the
connected components of $E \times_c \Gamma$ are mutually isomorphic by the
transitivity of the natural action of $\Gamma$ on $E \times_c \Gamma$
-- which is certainly not the case in general.

\begin{dfn} \label{Tvolt}
Let $E$ be a directed graph, $c : E^1 \rightarrow \Gamma$ a
function where $\Gamma$ is a countable group. Fix $T$ a spanning
tree of $E$ and vertex $v \in E^0$ then define the {\em
$T$--voltage based at $v$} to be the function $c_{v,T} : E^1
\rightarrow \Gamma$ given by
\begin{equation} \label{Tvoltdef}
c_{v,T} ( e ) := c ( b_{s(e)} e b_{r(e)}^{-1} )  = c ( b_{s(e)} ) c (e)
c ( b_{r(e)} )^{-1} ,
\end{equation}

\noindent
where for $w \in E^0$, $b_w$ is the unique reduced walk in $T$ from
$v$ to $w$.
\end{dfn}

\noindent The function $c_{v,T} : E^1 \rightarrow \Gamma$ is the
same as the one used in the proof of Theorem \ref{coverchar}. The
following result is not difficult to establish

\begin{prop} \label{Tvoltprop}
Let $E$ be a directed graph with spanning tree $T$ and $c : E^1 \rightarrow
\Gamma$ be a function where $\Gamma$ is a countable group. Then for $v
\in E^0$, we have $c_{v,T} \sim c$;  furthermore $\Gamma_v (c) =
\Gamma_v ( c_{v,T} )$.
\end{prop}

\begin{rmks} \label{indep}
Proposition \ref{Tvoltprop} shows that if we use
another base vertex to compute the $T$-voltages or if we use a
different spanning tree then we get an cohomologous function from
$E^1$ to $\Gamma$. Moreover, if $T$ is a spanning tree
for $E$ then the local voltage group of
$c$ at $v$ is generated by $\{ c_{v,T} (e) : e \in E^1 \backslash T^1
\}$ (cf.\ \cite[Theorem 2.5.3]{gt}).
\end{rmks}

\begin{prop} \label{compspg}
Let $E$ be a connected directed graph with spanning tree $T$, $v \in E^0$
and $c : E^1 \rightarrow \Gamma$ a function where $\Gamma$ is a countable
group. Then each connected component of
$E \times_c \Gamma$ is isomorphic to $E \times_{c_{v,T}} \Gamma_v (c)$.
\end{prop}

\begin{proof}
For $w \in E^0$ let $b_w$ be the unique reduced walk in $T$ from
$v$ to $w$. Since each of the connected components of $E \times_c
\Gamma$ are isomorphic by Corollary \ref{my2-2}, we shall deal with the
connected component $F$ which contains the vertex $(v, 1_\Gamma
)$. Define $\phi : E \times_{c_{v,T}} \Gamma_v (c) \rightarrow E
\times_c \Gamma$ by 
\[
\phi^0 (v,c(a)) = ( v, c (a) c ( b_v ) ) \text{ and }
\phi^1 (e,c(a)) = ( e , c (a) c ( b_{s(e)} ) ) ,
\]

\noindent
where $a \in \pi_1 (E,v)$. We check that $\phi$ is a graph morphism:
\begin{align*}
\phi^0 ( r (e,c(a) ) = \phi^0 ( r(e) c(a) c ( b_{s(e)} e b_{r(e)}^{-1} ) )
&= ( r(e) , c(a) c ( b_{s(e)} e ) ) \\
&= r ( e , c(a) c ( b_{s(e)} ) = r ( \phi^1 ( e,c(a)) ,
\end{align*}

\noindent
and
\[
\phi^0 ( s(e,c(a))) = (s(e) , c(a) b_{s(e)} ) = s ( e , c(a) b_{s(e)}
) = s ( \phi^1 ( e,c(a) ) . 
\]

\noindent
To see that $\phi$ is injective, suppose that $\phi^1 (e,c(a)) =
\phi^1 (f,c(b))$, then $(e,c(a) c ( b_{s(e)} ) ) = (f, c(b) c (
b_{s(f)} ))$ in which case $e=f$ and then $c(a)=c(b)$. To show that
$\phi^0$ is injective is similar. It only remains to show that the
image of $\phi$ is $F$. Suppose that $(u,g) \in F^0$ then
there is a reduced walk $\tilde{b}$ from $(v, 1_\Gamma )$ such that $b
= p_c ( \tilde{b} )$ is a reduced walk in $E$ from $v$ to $u$ with $g =
c(b)$. Since $b b_u^{-1}$ is a closed reduced walk in $E$ which begins at $v$
it follows that $c ( b b_u^{-1} ) \in \Gamma_v (c)$ and then
$(u, c(b b_u^{-1} )$ is a vertex in $E \times_{c_{v,T}} \Gamma_v (c)$
such that $\phi^0 ( u, c ( b b_u^{-1} ) ) = (u,c(b))=(u,g)$. A similar
proof shows that $\phi^1$ is surjective which completes our proof.
\end{proof}

\begin{rmk} \label{covercofinal}
Let $E$ be a directed graph and $c : E^1 \rightarrow \Gamma$ a
function. If $E \times_c \Gamma$ is cofinal then $E$ is cofinal,
however the converse is not true.
\end{rmk}

\noindent We shall be interested in the case when each edge in $E$
is labelled by the integer $1$. The skew product graph $E \times_c
{\bf Z}$ then has no loops (see \cite[Proposition 2.6]{kp}). The
next result can be deduced from the definitions and the results
from section \ref{nnm}.

\begin{prop} \label{skewcong}
Let $E$ be strongly connected  and let $c : E^1 \rightarrow {\bf Z}$
given by $c(e) = 1$ for all $e \in E^1$. If $E$ has period $d$ then
$\Gamma_v ( c ) \cong d {\bf Z}$ for all $v \in E^0$.
\end{prop}

\noindent If $E^0$ is finite then we can say a little more.

\begin{prop} \label{compcofinal}
Let $E$ be a strongly connected graph of period $d$ with finitely
many vertices. Let $c : E^1 \rightarrow {\bf Z}$ be given by $c
(e) = 1$ for all $e \in E^1$, then each connected component of $E
\times_c {\bf Z}$ is cofinal.
\end{prop}

\begin{proof}
Suppose that $\tilde{x}$ is an infinite path in one connected
component of $E \times_c {\bf Z}$ and $(v,m)$ a vertex in the same
component. Since $E^0$ is finite there is $u \in E^0$ and an
increasing sequence $n_1 \ge n_2 \ge n_3 \ge \ldots$ of integers
such that $( u , n_i )$ is the range of some edge in $\tilde{x}$.
Let $\alpha$ be a shortest path from $v = p_c ( v , m )$ to $w = s
( p_c ( \tilde{x}) )$ in $E$, then $\vert \alpha \vert = r_w$, the
residue class of $w$ in $E$ with respect to $v$ (see Lemma
\ref{resclass}). Let $\tilde{\alpha}$ be the lift of $\alpha$ in
$E \times_c {\bf Z}$ with source $(v,m)$, then by Remark
\ref{walks}, $r ( \tilde{\alpha} ) = ( w , m + r_w )$. Since
$(w,n) := s ( \tilde{x} )$ and $(w , m + r_w )$ belong to the same
connected component, by Proposition \ref{logrcor} we must have
\begin{equation} \label{link1}
n -  ( m + r_w ) = k' d
\end{equation}

\noindent
for some $k' \in {\bf Z}$. Since $(w,n)$ connects to $( u , n_i )$, it
also follows from Lemma \ref{resclass}(ii) that for all $i \ge 1$ we have
\begin{equation} \label{link2}
n_i = n + k_i d + s_u
\end{equation}

\noindent
for some $k_i \ge 1$, where $s_u$ is the residue class of $u$ in $E$
with respect to $w$. By Lemma \ref{resclass}(ii) there is a positive
integer $N (u)$ such that for all $k \ge N(u)$ there exists a path of
length $k d + s_u$ from $w$ to $u$. Since the $k_i$'s become
arbitrarily large we may assume that $k'+k_i \ge N ( u )$ for
sufficiently large $i$, and so by Lemma \ref{resclass}(ii)  there
exists paths $\beta_i$ in $E$ from $w$ to $u$ of length $(k'+k_i)d +
s_u$. Let $\tilde{\beta}_i$ be the lift of $\beta_i$ in $E \times_c {\bf
Z}$ with source $( w , m + r_w )$, then by Remark \ref{walks}
\begin{align*}
r ( \tilde{\beta}_i ) &= (  u , m + r_w + (k'+k_i)d + s_u ) \\
&= (  u , m + r_w + ( n - ( m + r_w ) ) + k_i d + s_u ) \text{ by (}
\ref{link1}
\text{)} \\
&= ( u , n + k_i d + s_u ) = ( u , n_i ) \text{ by (} \ref{link2} \text{).}
\end{align*}

\noindent
Hence $\tilde{\alpha} \tilde{\beta}_i$ is a path in the connected
component of $E \times_c {\bf Z}$ containing $(v,m)$ and $(w,n)$
beginning at $(v,m)$ and ending at $( u ,  n_i )$ for sufficiently
large $i$. Since $\tilde{x}$ passes through $( u , n_i )$ for all $i$
it follows that this component is cofinal by definition.
\end{proof}

\section{The AF core}

\noindent
Let $E$ be a directed graph, and $\{ s_e , p_v \}$ a Cuntz-Krieger
$E$-family generating $C^* (E)$. There is a strongly continuous
action $\gamma$ of ${\bf T}$ on $C^* (E)$ called the {\em gauge
action} which is characterised by
\[
\gamma_z s_e = z s_e , \text{ and } \gamma_z p_v = p_v ,
\]
where $z \in {\bf T}$. The fixed point algebra of the gauge
action, $C^* (E)^\gamma$ is AF and normally referred to as the AF
core (see \cite[\S 2.2]{pr} and \cite[\S 2]{bprs}). We wish to
establish conditions which guarantee that $C^* (E)^\gamma$ is
simple. To do this we shall use the following Morita equivalence which
was established in \cite[Proposition 2.8]{kp}.

\begin{lem} \label{fpaandspg}
Let $E$ be a row-finite directed graph with no sinks and sources, and
$c : E^1 
\rightarrow {\bf Z}$ be given by $c (e) = 1$ for all $e \in E^1$.
Then $C^* (E)^\gamma$ is strongly Morita equivalent to $C^* ( E
\times_c {\bf Z} )$.
\end{lem}

\noindent
The case when there are finitely many vertices is straightforward.

\begin{thm} \label{fpafinitecase}
Let $E$ be a row-finite directed graph with no sinks or sources and $E^0$
finite. Then $C^* (E)^\gamma$ is simple if
and only if $E$ is strongly connected with period $1$.
\end{thm}

\begin{proof}
Let $E$ be strongly connected with period one, and $c : E^1
\rightarrow {\bf Z}$ be given by $c(e)=1$ for all $e \in E^1$. By
Proposition \ref{compcofinal} the skew product graph $E \times_c {\bf Z}$ is
cofinal. Since $E \times_c {\bf Z}$ has no loops and $E$ is row-finite, by
Theorem \ref{ps} $C^* ( E \times_c {\bf Z} )$ is simple, and
then so is $C^* ( E )^\gamma$ by Lemma \ref{fpaandspg}.

Now suppose that $C^* (E)^\gamma$ is simple, then by Lemma
\ref{fpaandspg} $C^* ( E \times_c {\bf Z} )$ is simple. From Theorem
\ref{ps} $E \times_c {\bf Z}$ is cofinal and so $E$ is cofinal by
Remark \ref{covercofinal}. We claim that $E$ is strongly
connected. Let $u , v \in E^0$ then since $E^0$ is finite and $E$ has
no sources there is a loop $\alpha \in E^*$ and a path $\beta \in E^*$
with $s ( \alpha ) = s ( \beta )$ and $r ( \beta ) = v$. Let $x =
\alpha \alpha \cdots \in E^\infty$ then since $E$ is
cofinal there is a path $\mu \in E^*$ with $s ( \mu ) = u$ and $r (
\mu ) =  s ( \alpha )$. Hence there is a path $\mu \beta \in E^*$ from
$u$ to $v$, which establishes our claim. From Proposition \ref{skewcong} and
Corollary \ref{my2-2} it follows that $E$ has period $1$.
\end{proof}

\begin{rmk}
Theorem \ref{fpafinitecase} is not true if $E^0$ is infinite. The
graph $E$ shown below
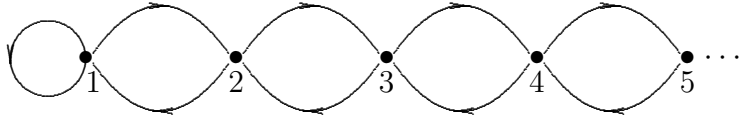
\[
\beginpicture
\setcoordinatesystem units <1cm,0.7cm>
\setplotarea x from -0.5 to 8.5, y from -1 to 1

\circulararc 360 degrees from 0 0 center at -0.5 0
\arrow <0.25cm> [0.1,0.3] from -1 0  to -1 -0.1

\setquadratic
\plot 0.1 0.1 1 1 1.9 0.1 /
\plot 0.1 -0.1 1 -1 1.9 -0.1 /

\arrow <0.25cm> [0.1,0.3] from 1 1  to 1.1 1
\arrow <0.25cm> [0.1,0.3] from 1 -1 to 0.9 -1

\setquadratic
\plot 2.1 0.1 3 1 3.9 0.1 /
\plot 2.1 -0.1 3 -1 3.9 -0.1 /

\arrow <0.25cm> [0.1,0.3] from 3 1  to 3.1 1
\arrow <0.25cm> [0.1,0.3] from 3 -1 to 2.9 -1

\setquadratic
\plot 4.1 0.1 5 1 5.9 0.1 /
\plot 4.1 -0.1 5 -1 5.9 -0.1 /

\arrow <0.25cm> [0.1,0.3] from 5 1  to 5.1 1
\arrow <0.25cm> [0.1,0.3] from 5 -1 to 4.9 -1

\setquadratic
\plot 6.1 0.1 7 1 7.9 0.1 /
\plot 6.1 -0.1 7 -1 7.9 -0.1 /

\arrow <0.25cm> [0.1,0.3] from 7 1  to 7.1 1
\arrow <0.25cm> [0.1,0.3] from 7 -1 to 6.9 -1

\put{$\cdots$} at 8.5 0
\put{$\bullet$} at 0 0
\put{$\bullet$} at 2 0
\put{$\bullet$} at 4 0
\put{$\bullet$} at 6 0
\put{$\bullet$} at 8 0

\put{$1$}[t] at 0.1 -0.25
\put{$2$}[t] at 2 -0.25
\put{$3$}[t] at 4 -0.25
\put{$4$}[t] at 6 -0.25
\put{$5$}[t] at 8 -0.25

\endpicture
\]

\noindent
is clearly strongly connected with period $1$. However
for the function $c : E^1 \rightarrow {\bf Z}$ given by $c(e)=1$ for
all $e \in E^1$, the skew product graph $E \times_c {\bf Z}$ is not
cofinal.
\end{rmk}

\noindent
If $E$ has sinks or a vertex of infinite valency then we
will make use of the following result:

\begin{lem} \label{gaugecase}
Let $E$ be a directed graph, and suppose that $H \subset E^0$ is a
nonempty saturated hereditary subset of $E^0$ such that $H \neq E^0$,
then $C^* (E)^\gamma$ is not simple.
\end{lem}

\begin{proof}
Let $H$ be a nontrivial saturated hereditary subset of $E^0$ and $I_H$
be the gauge invariant ideal of $C^* (E)$ generated by $\{
p_v : v \in H \}$. Let $\Phi : C^* (E) \rightarrow C^* ( E )^\gamma$
be the conditional expectation associated to the gauge action, then
$\Phi$ is faithful on positive elements (cf.\ \cite[\S 1]{bprs}). Since $H$
is a proper subset of $E^0$, and $\Phi$ an expectation which is
faithful on positive elements, it follows that
$\Phi ( I_H )$ is a  proper ideal of $C^* (E)^\gamma$ and the result follows.
\end{proof}

\noindent
Let us first deal with the case when there are sinks.

\begin{prop} \label{noloopscase}
Let $E$ be a row-finite directed graph which has at
least one sink. Then $C^* (E)^\gamma$ is simple if and only if $E$
consists of a single vertex.
\end{prop}

\begin{proof}
If $E$ consists a single vertex, then since
are no edges we have $C^* ( E) = C^* (E)^\gamma = {\bf C}$ which is simple.
Suppose that $E$ has sinks $v_1 \neq v_2$ then $\Sigma  ( \{ v_1 \} )$
is a proper saturated hereditary subset of $E^0$ and so gives rise to
a proper ideal in $C^* (E)^\gamma$ by Lemma \ref{gaugecase}.
Suppose now that $v$ is the only sink and let $J$ be the subalgebra of
$C^* (E)^\gamma$
generated by elements of the form $s_\mu s_\nu^*$ where $\mu , \nu \in E^n$
for some $n \ge 1$ are such that $r ( \mu ) = r ( \nu )$. It is
straightforward 
to show that $J$ is a closed $2$-sided ideal in $C^* (E)$.
If $E$ has edges then $J$ is nontrivial and $J \neq C^* (E)$
since $p_v \not\in J$. Hence $C^* (E)^\gamma$
is simple only if there are no edges and $E^0 = \{ v \}$.
\end{proof}

\noindent
The case when there is a vertex of infinite valency is quite similar.

\begin{prop} \label{infvertcase}
Let $E$ be a directed graph with at least one
vertex of infinite valency, then $C^* (E)^\gamma$ is not simple.
\end{prop}

\begin{proof}
Suppose $v \in E^0$ has infinite valency. Let
$J$ be the $C^*$-subalgebra of $C^* (E)^\gamma$ generated by elements
of the form $s_\mu s_\nu^*$ where $\mu , \nu \in E^n$ for some $n \ge 1$
are such that $r ( \mu ) = r ( \nu )$. It is straightforward to show
that $J$ is a closed $2$-sided ideal in $C^*
(E)^\gamma$. Now $J \neq C^* (E)^\gamma$ since $p_v \not\in J$, and
$J$ is nonzero since $E$ has edges, so the result follows.
\end{proof}

\noindent
Now suppose that $E^0$ is infinite and $E$ has no sinks.
If $E$ is either strongly connected or has no loops
then we claim that for each $v \in E^0$ there is an
infinite path $x(v)$, with source $v$ such
that for all $n$ there is no path of length less than $n$ from $v$ to
$r ( x (v)_n ) $. Every vertex on such a path is visited in the
shortest distance from $v$. An infinite path with this property is
said to be $v$-{\em depth-first}, since it is a path in the
depth-first spanning tree of $E$ (see \cite[\S 3.3]{co}).

For $v \in E^0$ let $V (0) = \{ v \}$ and for $n \ge 1$ let $V(n)$
denote those vertices to which $v$ connects by a path whose length is
less than $n$. Since $E^0$ is infinite with no sinks and $E$ is strongly 
connected or has no loops it follows that $V(n)$ is nonempty and
the containment $V (n-1) \subset V(n)$ is strict for $n \ge
1$. The set $V(n) \backslash V(n-1)$ denotes those
vertices which can be reached from $v$ with shortest path of length $n$.
We say that a path $\alpha \in E^n (v) := \{ \alpha \in
E^n : s ( \alpha ) = v \}$ is $v$-{\em deep} if it uses different vertices,
that is if $r ( \alpha_i ) \in V(i) \backslash V(i-1)$ for $1 \le i
\le n$. Let $E^n_g  (v)$ denote the collection of $v$-deep paths in $E^n (v)$.

\begin{lem} \label{indstep}
Let $E$ be a directed graph with $E^0$ infinite and no sinks.
If $E$ is strongly connected or has no loops then for $v \in E^0$
the set $E^n_g (v)$ is non-empty for all $n \ge 1$.
\end{lem}

\begin{proof}
The result then follows immediately from the fact that $r \left( E^n_g (v)
\right) = V(n) \backslash V(n-1)$ which is non-empty for all $n \ge 1$.
\end{proof}

\begin{thm}
Let $E$ be a row-finite directed graph with $E^0$ infinite and no
sinks. If $E$ is strongly connected or has no loops then for each
$v \in E^0$ there is a $v$-depth-first path.
\end{thm}

\begin{proof}
For integers $m > n \ge 0$ let $E^{n,m}_g (v)$ denote those $\alpha'
\in E^n_g (v)$ which appear as the first part of some $\alpha \in
E^m_g (v)$. Clearly $E^{n,m}_g (v) \neq \emptyset$ for all $m > n \ge
0$ since $E^m_g (v)$ is non-empty and for any $\alpha = \alpha'
\alpha'' \in E^m_g (v)$ with $\alpha' \in E^n (v)$ we have
$\alpha' \in E^{n,m}_g (v)$. Moreover, this argument also shows that
$E^{n,m}_g (v) \supseteq E^{n,p}_g (v)$ for all integers $p > m > n
\ge 0$. Let
\[
E^{n,\infty}_g (v) = \bigcap_{m > n} E^{n,m}_g (v) ,
\]

\noindent
then $E^{n,\infty}_g (v) \neq \emptyset$ since it is the intersection
of a decreasing sequence of finite non-empty sets. Now we may define the
infinite path we seek: Since $E^{1,\infty}_g (v) \neq \emptyset$ there
is $x(v)_1 \in E^1$ such that there are $\alpha_n \in E^n$ for $n \ge
1$ with  $x(v)_1 \alpha_n \in E^{1,n+1}_g (v)$.
Since $E$ is row-finite there is $x(v)_2 \in E^1$ such that there are
infinitely many such $\alpha_n$ with $\alpha_n = x(v)_2 \alpha''$
where $\alpha'' \in E^{n-1}$. In particular $x(v)_1 x(v)_2 \in
E^{2,\infty}_g (v)$ since for $n > 2$ there are $\alpha \in E^{2,n}_g
(v)$ with $\alpha  = x(v)_1 x(v)_2 \alpha''$ for some $\alpha''  \in
E^{n-2}$. In this way we may define $x(v) = ( x(v)_i ) \in E^\infty$
with $s( x(v)_1 ) = v$ and $r ( x(v)_i ) \in V (i) \backslash V (i-1)$
for all $i \ge 1$.
\end{proof}

\begin{cor} \label{infcase}
Let $E$ be a row-finite directed graph with $E^0$ infinite and no
sinks. If $E$ is strongly connected or has no loops then $C^*
(E)^\gamma$ is not simple.
\end{cor}

\begin{proof}
Let $c : E^1 \rightarrow {\bf Z}$ be given by $c(e)=1$
for all $e \in E^1$ and for $v \in E^0$, let $x(v)$ be a
$v$-depth-first path. Let $x' \in ( E \times_c {\bf Z} )^\infty$ be the
lift of $x(v)$ with source $(-1,v)$, and suppose that $\alpha'$ is such
that $s ( \alpha ' )  =  ( 0 , v )$ and $r ( \alpha' ) = r ( {x'}_i ) = (
i-1, r ( x(v)_{i} ) )$ for some $i$. Hence  $\alpha = p_c ( \alpha' )$ is a
path of length $i-1$  from $v$ to $ r ( x (v)_{i} )$ which means that
$r ( x (v)_i ) \in V(i-1)$, which contradictions the definition
of $x(v)$. Hence $E \times_c {\bf Z}$ is not cofinal,
so $C^* (E \times_c {\bf Z} )$ is not simple by Theorem \ref{ps}
and hence $C^* (E)^\gamma$
is not simple by Lemma \ref{fpaandspg}
\end{proof}

\begin{thm} \label{above}
Let $E$ be a directed graph, then $C^* (E)^\gamma$ is
simple if and only if either $E$ has a cofinal strongly connected subgraph
of period $1$ which has finitely many vertices or $E$ consists of a
single vertex.
\end{thm}

\begin{proof}
If $E$ has a vertex infinite valency then $C^* (E)^\gamma$ is not
simple by Proposition \ref{infvertcase}. So we may suppose that
$E$ is row-finite. If $E$ is not cofinal then, by the proof of
Theorem \ref{ps}, there is a nontrivial saturated hereditary
subset of $E^0$ which gives rise to a nontrivial ideal in $C^*
(E)^\gamma$ by Lemma \ref{gaugecase}. So we may suppose that $E$
is cofinal. If $E$ does not satisfy condition (K)
then there is a cofinal subgraph $L$ which consists of a loop with no
exits. Let $c : E^1 \rightarrow {\bf Z}$ be given by $c(e)=1$ for
all $e \in E^1$, then $E \times_c {\bf Z}$ is cofinal if and only
if $L$ has one vertex. Hence by Lemma
\ref{fpaandspg}, $C^* (E)^\gamma$ is simple if and only if
$E$ has a cofinal strongly connected subgraph of period $1$
which consists of a single vertex (and edge). So we may suppose that
$E$ also satisfies condition
(K). Therefore we may assume that $E$ satisfies conditions (i)-(iii)
of Theorem \ref{ps} and so $C^* (E)$ is simple.

If $E$ has a sink then by Proposition \ref{noloopscase},
$C^* (E)^\gamma$ is simple if and only if $E$
consists of a single vertex. So we suppose that $E$ has no sinks.
By Theorem \ref{reducetosc} either $E$ has no loops or there is a
strongly connected cofinal subgraph $F \subseteq E$. Suppose that
$E$ has no loops, then $E^0$ cannot be finite (as that would mean there
had to be a sink cf.\ \cite[\S 2]{kpr}).
If $E^0$ is infinite then $C^* (E)^\gamma$ is not simple by Corollary
\ref{infcase}. Suppose $E$ has a strongly connected cofinal subgraph $F$
then by Corollary \ref{ga-appl} there is a projection $P \in
M (C^*(E) )$ such that $C^* (F) \cong P C^* (E) P$. Since this
isomorphism is essentially the identity map, it
commutes with the usual gauge actions on $C^* (E)$ and $C^*
(F)$ which are both denoted by $\gamma$. The projection $P$ is the
limit of a sum
of projections in $C^* (E)$ which are invariant under the gauge
action. Hence $P \in M (C^* (E) )^\gamma$, and then $C^* (F)^\gamma
\cong P C^* (E)^\gamma P$. Since $F$ is cofinal in $E$, it follows that
$P C^* (E)^\gamma P$ is full and so $C^* (E)^\gamma$ is strongly Morita
equivalent to $C^* (F)^\gamma$. The result now follows from Theorem
\ref{fpafinitecase} and Theorem \ref{infcase}.
\end{proof}

\noindent
If we discount the trivial case when the graph only consists of a single
vertex, Theorem \ref{above} says the following: Up to Morita equivalence, the
only graphs $E$ for which the AF core of $C^* (E)$ is simple are
those which have finitely many vertices and are strongly
connected with period one.
For strongly connected graphs we have
the following structure results for the AF core.

\begin{thm}
Let $E$ be a strongly connected row-finite graph with period $d$. Then $C^*
(E)^\gamma$ is a direct sum of $d$ mutually isomorphic AF
algebras. If, in addition $E^0$ is finite then these AF algebras are simple.
\end{thm}

\begin{proof}
Let $E$ be strongly connected, row-finite with period $d$, and $c : E^1
\rightarrow {\bf Z}$ be given by $c (e) =1$ for all $e \in E^1$.
By Corollary \ref{my2-2} $E \times_c {\bf Z}$
consists of $d$ mutually isomorphic components. Since each component is a
subgraph of a graph which has no loops, the first part follows from
Lemma \ref{fpaandspg} and \cite[Theorem 2.4]{kpr}. If $E^0$ is finite
then each component is cofinal by Proposition \ref{compcofinal} and
the last statement follows from Theorem \ref{ps}.
\end{proof}

\begin{thm}
Let $E$ be a strongly connected row-finite graph with finitely many
vertices, then $C^* (E)$
is stably isomorphic to a crossed product of a simple AF algebra by an
action of ${\bf Z}$.
\end{thm}

\begin{proof}
Let $E$ be a strongly connected graph with period $d$, $T$ be a
spanning tree for $E$, $v \in E^0$ and $c : E^1 \rightarrow {\bf Z}$
be given by $c(e)=1$ for all $e \in E^1$. By Proposition
\ref{compspg} each component of $E \times_{c} {\bf Z}$ is
isomorphic to $E \times_{c_{v,T}} \Gamma_v (c )$. Since $E$ is
strongly connected with period $d$, by Proposition
\ref{skewcong} we
have $\Gamma_v (c) = d {\bf Z} \cong {\bf Z}$. Let
$\lambda$ denote the free ${\bf Z}$-action on $E \times_{c_{v,T}} d
{\bf Z}$ which has quotient $E$ then by \cite[Corollary 3.9]{kp} we have
\[
C^* ( E \times_{c_{v,T}} d {\bf Z} ) \times_\lambda {\bf Z} \cong C^*
(E) \otimes \mathcal{K} \big( \ell^2 ( {\bf Z} ) \big) .
\]

\noindent
Since $E \times_{c_{v,T}} d {\bf Z}$ is isomorphic to a subgraph of $E
\times_c {\bf Z}$ it has no loops, and so $C^* ( E \times_{c_{v,T}} d {\bf Z}
)$ is AF by \cite[Theorem 2.4]{kpr}. If $E^0$ is finite then $E
\times_{c_{v,T}} d {\bf Z}$ is  cofinal by Proposition
\ref{compcofinal}. It then follows that $C^* ( E
\times_{c_{v,T}} d  {\bf Z} )$ is  simple by Theorem \ref{ps}
which completes the proof.
\end{proof}

\end{document}